\documentclass[11pt,reqno,oneside]{amsart}
 \usepackage{latexsym}
 \usepackage{amssymb}
 \usepackage{amsfonts}
\usepackage{graphics}
 \author[Mathonet]{P. Mathonet}
\thanks{University of Li\`ege, Institute of mathematics, Grande Traverse, 12 - B37, B-4000 Li\`ege, Belgium.  email : P.Mathonet[at]ulg.ac.be}
\author[Radoux]{F. Radoux}
\thanks{Universit\'e du Luxembourg, URM, Campus Limpertsberg, Avenue de la fa\"{i}encerie 162 A, L-1511 Luxembourg. email : Fabian.Radoux[at]uni.lu}
\date{\today}

\newtheorem{lem}{Lemma}
\newtheorem{thm}[lem]{Theorem}
\newtheorem{prop}[lem]{Proposition}
\newtheorem{cor}[lem]{Corollary}
\theoremstyle{remark}
\newtheorem{rem}{Remark}
\theoremstyle{definition}
\newtheorem{defi}{Definition}

\newcommand{\R}{\mathbb{R}}

\renewcommand{\L}{\mathcal{L}}

\renewcommand{\S}[1][\delta]{\mathcal{S}_{{#1}}}

\newcommand{\g}{\mathfrak{g}}

\newcommand{\h}{\mathfrak{h}}

\newcommand{\cc}{{\mathcal{C}}}
\renewcommand{\L}{\mathcal L}
\newcommand{\om}{\omega^{-1}}
\renewcommand{\k}{\kappa}

\title{On Natural and Conformally Equivariant Quantizations}
\begin{document}
\begin{abstract}
The concept of conformally equivariant quantizations was introduced by Duval,
Lecomte and Ovsienko in \cite{DLO} for manifolds endowed with flat conformal 
 structures. They obtained results of existence and uniqueness (up to
 normalization) of such a quantization procedure. A natural generalization
 of this concept is to seek for a quantization procedure, over a manifold $M$, 
 that depends on a pseudo-Riemannian metric, is natural and is
 invariant with respect to a conformal change of the metric. The existence of such a procedure was conjectured by P. Lecomte in \cite{Leconj} and proved by C. Duval and V. Ovsienko in \cite{DO1} for symbols of degree at most 2 and by S. Loubon Djounga in \cite{Loubon} for symbols of degree 3.
In two recent papers \cite{MR,MR1}, we investigated the question of existence of projectively equivariant quantizations using the framework of Cartan connections. Here we will show how the formalism developed in these works adapts in order to deal with the conformally equivariant quantization for symbols of degree at most 3.
 This will allow us to easily recover the results of \cite{DO1} and \cite{Loubon}. We will then show how it can be modified in order to prove the existence of conformally equivariant quantizations for symbols of degree 4.
\end{abstract}
\maketitle
\section{Introduction}
It is common in the mathematical literature to think of a quantization map as a linear bijection from a space of classical observables to a space of differential operators acting on wave
functions (see \cite{Woo}).
In our setting, the observables (also called \emph{Symbols}) are smooth functions on the cotangent bundle $T^*M$ of a manifold $M$, that are polynomial along the fibres, while the differential operators act on half densities.

The concept of \emph{equivariant quantization} was introduced by Lecomte and Ovsienko in \cite{LO} and developed in a series of papers \cite{DLO, DO,BM, BHMP}. The idea is that, when a Lie group $G$ acts on $M$ by local diffeomorphisms, it is sometimes natural to require that the quantization map be equivariant with respect to the lifted actions of $G$ on symbols and on differential operators.

In \cite{LO} and \cite{DO}, the authors considered the case of the projective group
$PGL(m+1,\R)$ acting on the manifold $M=\R^m$ by linear fractional
transformations. This leads to the notion of projectively equivariant
quantization or its infinitesimal counterpart, 
the $sl(m+1,\R)$-equivariant quantization. 
In \cite{DLO}, the authors considered the group $SO(p+1,q+1)$ acting on the
space $\R^{p+q}$ or on a manifold endowed with a flat conformal structure.
 They extended the problem by considering the space $\mathcal{D}_{\lambda,\mu}$
of differential operators mapping $\lambda$-densities into $\mu$-densities and
a suitable space of symbols $\mathcal{S}_{\mu-\lambda}$.

In both situations, the results were the existence and uniqueness of 
equivariant quantization provided the shift value $\delta = \mu - \lambda$
does not belong to a set of critical values. These results settled the problem projectively- and conformally- equivariant quantizations in the framework of manifold endowed with flat structures.

Then the particular case of symbols of degree at most three was considered in detail in \cite{Bou1,Bou2} in the projective case and in \cite{DO1,Loubon} in the conformal case. There the authors showed that the equivariant quantization procedure can be expressed using a connection (the Levi Civita connection associated to a pseudo Riemannian metric in the conformal case), in such a way that it only depends on the projective class of the connection or on the conformal class of the metric.

In \cite{Leconj}, P. Lecomte conjectured the existence of a quantization 
procedure depending on a torsion-free linear connection (resp. on a
pseudo-Riemannian metric), that would be
 natural (in all arguments) and that would be invariant with respect to a projective (resp. conformal) change of connection (resp. metric).

 In the projective case, the existence of such a \emph{Natural and equivariant
   quantization procedure} was proved by M. Bordemann 
in \cite{Bor}, using the notion of Thomas-Whitehead connection associated to a projective class of connections.
His construction was later adapted by S. Hansoul (see \cite{Hansoul,sarah}) in order to deal with other types of differential operators.
 
In \cite{MR,MR1}, we analysed the existence problem for a natural and projectively
 equivariant quantization using the theory of projective Cartan connections. 
Among our motivations, one was to use the similarity between the theory of
 projective Cartan connections and the theory of conformal ones (we refer the
 reader to \cite{Kobabook} where both cases are presented concurrently) 
in order to deal with the conformal case. 

In these papers, we showed how the formulae that were obtained in the flat situation could be easily modified in order to obtain a natural and projectively equivariant quantization. 

In the present work, we will show that in the conformal situation, the same modifications allow to build a natural and conformally equivariant for symbols of degree at most three, thus recovering more easily the results of \cite{DO1,Loubon}. We will show that the same procedure of \cite{MR1} does not directly work for symbols of degree 4. We will nevertheless exhibit an explicit formula for the conformally equivariant quantization of such symbols, in terms of the Cartan connection associated to the conformal class of the metric. Moreover the rather simple expression of the solution might show the direction for a solution of the existence problem for symbols of arbitrary degree. 
\section{Basic notions and problem setting}
In this section, we will just fix some notation concerning tensor densities, symbols and differential operators. For more details about these notions, we refer the reader to \cite{MR,MR1} and references therein. Throughout this paper, we let $M$ denote an $m$-dimensionsal Hausdorff, second countable smooth manifold.
\subsection{Tensor densities}
Denote by $\Delta^{\lambda}(\R^m)$ the one dimensional representation of $GL(m,\R)$ given by
\[\rho(A) e = \vert det A\vert^{-\lambda} e,\quad\forall A\in
GL(m,\R),\;\forall e\in \Delta^{\lambda}(\R^m).\]
The vector bundle of $\lambda$-densities is then defined by
\[P^1M\times_{\rho}\Delta^{\lambda}(\R^m)\to M,\]
where $P^1M$ is the linear frame bundle of $M$.

Recall that the space $\mathcal{F}_{\lambda}(M)$ of smooth sections of this bundle, the space of $\lambda$-densities,  can be identified with the space $C^{\infty}(P^1M,
\Delta^{\lambda}(\R^m))_{GL(m,\R)}$ of functions $f$ such that
\[f(u A) = \rho (A^{-1}) f(u)\quad \forall u \in P^1M,\;\forall A\in
GL(m,\R).\]
Moreover, since the bundle of $\lambda$-densities is associated to the linear frame bundle, there are natural actions of local diffeomorphisms and of vector fields on ${\mathcal F}_{\lambda}$. These actions were detailed for instance in \cite{DLO,LO}.
\subsection{Differential operators and symbols}
As usual, we denote by $\mathcal{D}_{\lambda,\mu}(M)$ the space 
of differential operators from $\mathcal{F}_{\lambda}(M)$ to
$\mathcal{F}_{\mu}(M)$. The actions of vector fields and of (local) diffeomorphisms on $\mathcal{D}_{\lambda,\mu}(M)$ are induced by the actions on $\mathcal{F}_{\lambda}(M)$ and $\mathcal{F}_{\mu}(M)$.

The space $\mathcal{D}_{\lambda,\mu}$ is filtered by 
the order of differential operators. This filtration is preserved by the action of local diffeomorphisms and of vector fields.
The space of \emph{symbols} is then the associated graded space of $\mathcal{D}_{\lambda,\mu}$. It is also known that the principal operator $\sigma$ allows to identify the space of symbols with the space of contravariant symmetric tensor fields with coefficients in $\delta$-densities where $\delta=\mu-\lambda$ is the shift value.

More precisely, we denote by $S^l_{\delta}(\R^m)$ or simply $S^l_{\delta}$ the space $S^l\R^m\otimes\Delta^{\delta}(\R^m)$ endowed with the natural representation of $GL(m,\R)$. Then the vector bundle of symbols of degree $l$ is 
\[P^1M\times_{\rho}S^l_{\delta}(\R^m)\to M.\]
The space  $\mathcal{S}^l_{\delta}(M)$ of symbols of degree $l$ is then the space of smooth sections of this bundle, which can be identified with $C^{\infty}(P^1M,
S^l_{\delta}(\R^m))_{ GL(m,\R)}$.
Finally, the whole space of symbols is
\[\mathcal{S}_{\delta}(M)
=\bigoplus_{l=0}^{\infty}\mathcal{S}^l_{\delta}(M),\]
endowed with the classical actions of diffeomorphisms and of vector fields.
\subsection{Natural and equivariant quantizations}
A \emph{quantization on $M$} is a linear bijection $Q_{M}$
  from the space of symbols $\mathcal{S}_{\delta}(M)$ to the space of differential operators $\mathcal{D}_{\lambda,\mu}(M)$ such
  that
\[\sigma(Q_{M}(S))=S,\;\forall S\in\mathcal{S}_{\delta}^{k}(M),\;\forall
  k\in\mathbb{N},\]
where $\sigma$ is the principal symbol operator.

In the conformal sense, a \emph{natural quantization} is a collection of quantizations $Q_M$ depending on a pseudo-Riemannian metric  such that 
\begin{itemize}
\item
For all pseudo-Riemannian metric $g$ on $M$, $Q_{M}(g)$ is a quantization,
\item
If $\phi$ is a local diffeomorphism from $M$ to $N$, then one has
\[Q_{M}(\phi^{*}g)(\phi^{*}S)=\phi^{*}(Q_{N}(g)(S)),\]
 for all pseudo-Riemannian metrics $g$ on $N$, and all
$S\in\mathcal{S}_{\delta}(N).$
\end{itemize}
Recall now that two pseudo Riemannian metrics $g$ and $g'$ on a manifold $M$ are conformally equivalent if and only if there exists a positive function $f$ such that $g'=fg$. 

A quantization $Q_M$ is then \emph{conformally equivariant} if one has $Q_M(g)=Q_M(g')$  whenever $g$ and $g'$ are conformally equivalent .
\subsection{Conformal structures}\label{projconf}
In this section, we will recall the notions of conformal structures over a manifold $M$ and how they are associated to reductions of the second order frame bundle $P^2M$ to a certain group $H$. The description of this group as a semi direct product as well as the description of its Lie algebra is one of the main ingredients of our construction. This description was given for the projective and the conformal situation by S. Kobayashi in \cite{Kobabook}. Here we will recall the main results of \cite{Kobabook} and refer the reader to this book for more details. We will adopt the notation of \cite{MR,MR1} for the comparison with the projective case to be easy.

A conformal structure on a manifold is classically defined as an equivalence class of pseudo-Riemannian metrics. 
We now recall the group theoretic point of view.

Given $p$ and $q$ such that $p+q=m$, we consider the bilinear symmetric form of signature $(p+1,q+1)$ on $\R^{m+2}$ defined by
\[B:\R^{m+2}\times\R^{m+2}\to\R:(x,y)\mapsto {^ty}Sx,\]
where $S$ is the matrix of order $m+2$ given by
\[S=
      \left(
        \begin{array}{ccc}
           0 & 0 & -1\\
           0 & J & 0\\
          -1 & 0 & 0
       \end{array}
     \right).
    \]
Also the matrix
\[J=
      \left(
        \begin{array}{cc}
           I_{p} & 0\\
           0     & -I_{q}
          \end{array}
     \right)
    \]
    represents a nondegenerate symmetric bilinear form $g_0$ on $\R^m$, namely
    \[g_0:\R^{m}\times\R^{m}\to\R:(x,y)\mapsto {^ty}Jx.\]
    We also denote by $|x|^2$ the number $g_0(x,x)$. 
    
As we continue, we will use the classical isomorphism between $\R^m$ and $\R^{m*}$ defined by the symmetric bilinear form represented by $J$ :
\[\sharp : \R^m\to \R^{m*}: x\mapsto x^\sharp : x^\sharp(y)= g_0(x,y)\]
and $\flat=\sharp^{-1}$. 

The M\"{o}bius space is the projection of the light cone associated to this metric on the projective space $\mathbb{R}P^{m+1}$.
 
The group $G$ is made of linear transformations that leave $B$ invariant, modulo its center, that is,
\[G=\{X\in GL(m+2,\mathbb{R}) : {^tX}SX=S\}/\{\pm I_{m+2}\}.\]
It acts transitively on the M\"{o}bius space $S^{m}$.

The group $H$ is the
isotropy subgroup of $G$ at the point $[e_{m+2}]$ of the M\"{o}bius space :
\[H=\{ \left( 
           \begin{array}{ccc}
              a^{-1} & 0 & 0\\
              a^{-1}A\xi^\flat & A & 0\\
              \frac{1}{2a}|\xi|^2 & \xi & a
           \end{array}
        \right): A\in
        O(p,q),a\in\mathbb{R}_{0},\xi\in\mathbb{R}^{m*}\}/\{\pm I_{m+2}\}.\]
As in the projective situation, $H$ is a semi-direct product $G_{0}\rtimes
        G_{1}$. Here $G_{0}$ is isomorphic to $CO(p,q)$ and $G_{1}$ is
        isomorphic to $\mathbb{R}^{m*}$. There is also a projection 
\[\pi:H\mapsto CO(p,q):\left[
        \left( 
           \begin{array}{ccc}
              a^{-1} & 0 & 0\\
             a^{-1}A\xi^\flat & A & 0\\
               \frac{1}{2a}|\xi|^2 & \xi & a
           \end{array}
        \right)\right]\mapsto \frac {A}{a}.\]
The Lie algebra of $G$ is $\mathfrak{g}=so(p+1,q+1)$, and decomposes as a direct sum of subalgebras
\[\mathfrak{g}_{-1}\oplus\mathfrak{g}_{0}\oplus\mathfrak{g}_{1}\cong\mathbb{R}^{m}\oplus
co(p,q)\oplus\mathbb{R}^{m*}.\]
The isomorphism is given by
\[
         \left( 
           \begin{array}{ccc}
              -a & v^\sharp & 0\\
              \xi^{\flat} & A & v\\
              0 & \xi & a
           \end{array}
        \right)\mapsto(v,A-aI_{m},\xi).\]
This correspondence induces a structure of Lie algebra on $\mathbb{R}^{m}\oplus
co(p,q)\oplus\mathbb{R}^{m*}$. It is easy to see that the adjoint actions of $co(p,q)$ on $\R^m$ and on ${\R^m}^*$ coincides with the natural actions. Moreover, one has
\begin{equation}\label{bra} [h,x]=-x\otimes h + h^\flat\otimes x^\sharp -\langle h,x\rangle Id.\end{equation}

The Lie algebras corresponding to $G_{0}$,
$G_{1}$ and $H$ are respectively $\mathfrak{g}_{0}$, $\mathfrak{g}_{1}$, and
$\mathfrak{g}_{0}\oplus\mathfrak{g}_{1}$.

 It is well-known that there is a bijective and natural correspondence between
the conformal structures on $M$ and the reductions of $P^{1}M$ to the
structure group $G_{0}\cong CO(p,q)$. In \cite{Kobabook}, one shows that it is possible to
associate at each $G_{0}$-structure $P_{0}$ a principal $H$-bundle $P$ on $M$,
this association being natural and conformally invariant. Since $H$ can be considered as a subgroup of $G^2_m$, this $H$-bundle can
be considered as a reduction of $P^{2}M$. The relationship between conformal structures and reductions of $P^2M$ to $H$ is given by the following proposition.
\begin{prop}
There is a natural one-to-one correspondence between the conformal equivalence
classes of pseudo-Riemannian metrics on $M$ and the reductions of $P^2M$ to $H$.\end{prop}
Throughout this work, we will freely identify conformal structures and reductions of $P^2M$.
\subsection{Lift of equivariant functions}\label{Lift}
In the previous section, we recalled how to associate an $H$-principal bundle $P$ to a conformal structure $P_0$. We now recall how the densities and symbols can be regarded as equivariant functions on $P$.

If $(V,\rho)$ is a representation of $G_0$, then we may extend it to a representation $(V,\rho')$ of $H$ by
\[\rho'=\rho\circ\pi.\]
Now, using the representation $\rho'$, we can recall the relationship between
equivariant functions on $P_{0}$ and equivariant functions on $P$ (see \cite{Capcart}): if we
denote by $p$ the projection $P\to P_{0}$ , we have
\begin{prop}
If $(V,\rho)$ is a representation of $G_0$, then the map
$$p^{*}:C^{\infty}(P_{0},V)\mapsto C^{\infty}(P,V):f\mapsto f\circ p$$
defines a bijection from $C^{\infty}(P_{0},V)_{G_0}$ to
$C^{\infty}(P,V)_{H}$.
\end{prop}
\subsection{Cartan connections}
The main ingredient in our construction of natural and equivariant quantization is the normal Cartan connection associated to a conformal structure. In this section we recall the definitions and results about Cartan conformal connections that we will use as we continue. Let us begin with the general definition :

Let $L$ be a Lie group and $L_{0}$ a closed subgroup. Denote by $\mathfrak{l}$
and $\mathfrak{l}_{0}$ the corresponding Lie algebras. Let $N\to M$ be a
principal $L_{0}$-bundle over $M$, such that $\dim M=\dim L/L_{0}$. A Cartan
connection on $N$ is an $\mathfrak{l}$-valued one-form $\omega$ on $N$ such
that
\begin{enumerate}
\item
If $R_{a}$ denotes the right action of $a\in L_0$ on $N$, then
$R_{a}^{*}\omega=Ad(a^{-1})\omega$,
\item
If $k^{*}$ is the vertical vector field associated to $k\in\mathfrak{l}_0$,
then $\omega(k^{*})=k,$
\item
$\forall u\in N$, $\omega_{u}:T_{u}N\mapsto\mathfrak{l}$ is a linear
  bijection.
\end{enumerate}
When considering in this definition the principal $H$- bundle $P$ of the previous section, and taking as group $L$ the group $G$ and for $L_0$ the group $H$, we obtain the definition of Cartan conformal connections.

We will need some properties of the curvature of Cartan connections. 

If $\omega$ is a Cartan connection defined on an $H$-principal bundle $P$, then its
curvature $\Omega$ is defined by
 \begin{equation}\label{curv} 
\Omega = d\omega+\frac{1}{2}[\omega,\omega].
\end{equation}
The connection $\omega$ is constructed (in \cite{Kobabook}) in such a way that the curvature has values in $\h=\g_0\oplus\g_1$. 

Moreover, it is known that the curvature vanishes on vertical vector fields, that is, vector fields of the form $h^*$ for $h$ in $\h$.

It is common to use the function $\kappa\in C^{\infty}(P,\wedge^2\g^*\otimes\g)$ defined as follows
\[\kappa(x,y)=\Omega(\om(x),\om(y)).\]
Since $\Omega$ takes values in $\g_0\oplus\g_1$, we can decompose 
$\kappa$ as $\kappa_0+\kappa_1$ in an obvious manner.

Finally, since the function $\k$ vanishes on $\g_0\oplus\g_1$, the functions $\k_0$ and $\k_1$ actually belong to $C^{\infty}(P,\wedge^2\g_{-1}^*\otimes\g_0)$ and $C^{\infty}(P,\wedge^2\g_{-1}^*\otimes\g_1)$ respectively. Also notice that, since we identified $\g_0$ to $co(p,q)$, the function $\kappa_0$ also belongs to $C^{\infty}(P,\wedge^2\g_{-1}^*\otimes gl(\g_{-1}))$ while $\kappa_1$ belongs to $C^{\infty}(P,\wedge^2\g_{-1}^*\otimes\g_{-1}^*)$, since $g_1\cong {\R^m}^*$.

Moreover, if we denote by $(e_1,\cdots,e_m)$ the standard basis of $\g_{-1}\cong \R^m$, the components of the function $\kappa_0$ in this basis can be defined by
\[\k_{ijk}^l=(\kappa_0(e_j,e_k))_i^l.\]
Then a \emph{normal Cartan connection} is a Cartan connection whose curvature fulfills the relation
\[\k_{jil}^i=0\quad\forall j,l.\]
The following result (see \cite{Kobabook,Capinv}) ensures the existence of normal Cartan connections  
\begin{prop}
A unique normal conformal Cartan connection with values in the algebra $\mathfrak{g}$ is associated to every conformal structure $P$. This association is natural.
\end{prop}

We end this section by technical results about Lie derivatives of the curvature.
The following lemmas follow easily from the definition of the curvature and from the Ad-invariance of $\omega$.
\begin{lem}\label{omega}
There holds
\[[\om(x),\om(A)]=-\om(\k(x,A))+\om(L_{\om(x)}A)+\om([x,A]),\]
for every $x\in\g$ and $A\in C^{\infty}(P,\g).$

In particular there holds
\[ [h^*,\om(y)]=\om([h,y]),\]
for every $h\in\h$ and $y\in\g$.
\end{lem}
Finally we have
\begin{lem}\label{derkappa}
For every $h\in\g_0\oplus\g_1$, one has
\[L_{h^*}\omega=-ad(h)\circ\omega\quad\mbox{and}\quad L_{h^*}\Omega=-ad(h)\circ\Omega.\]
The function $\k$ is $\g_0$-invariant : for every $h\in\g_0$ one has
\[L_{h^*}\k=-ad(h)\k,\]
that is,
\[L_{h^*}\k(x,y)=-ad(h)(\k(x,y))+\k(ad(h)x,y)+\k(x,ad(h)y).\]
Finally, for every $h\in\g_1$, one has 
\[L_{h^*}\k_0=0\quad\mbox{and}\quad L_{h^*}\k_1(x,y)=\,[\k_0(x,y),h].\]
\end{lem}
\section{Equivariant quantizations in the flat case}
Our construction of the natural and conformally equivariant quantization is based on an adaptation of two previous techniques : the first one is the conformally equivariant quantization in the flat case (see \cite{DLO}) and the second is our construction of the projectively equivariant natural quantization in general. Here we recall the main features of the construction of the quantization in the flat case. The details and the proofs can be found in \cite{DLO,BM,MR1}.
\subsection{The flat case} 
Let us first recall the problem of equivariant quantizations over the Euclidean space $\R^m$.

In this framework, we make the following identifications :
\[\begin{array}{ccc}\mathcal{F}_{\lambda}(\R^m)&\cong&C^{\infty}(\R^m,\Delta^{\lambda}(\R^m)),\\
\mathcal{S}^k_{\delta}(\R^m) & \cong & C^{\infty}(\R^m,S^k_{\delta}(\R^m)).
\end{array}\]
The Lie algebra $\mathit{Vect}(\R^m)$ acts on these spaces in a well-known
manner : for every $X\in\mathit{Vect}(\R^m)$ and any symbol $s$, one has
\[L_Xs(x)=X.s(x) - \rho_{*}(D_xX)s(x)\]
where $\rho$ is the natural action of $GL(m,\R)$ on the fibre.
The space $\mathcal{D}_{\lambda \mu}(\R^m)$ of differential operators 
is equipped with the Lie derivative given by the commutator.
\subsubsection{Conformal algebra of vector fields}
The data of sections \ref{projconf} allow to define Lie algebras of vector fields. The Lie group $G$ acts on the homogeneous space $G/H$. Since this space can locally be identified with $\R^m$, there is a local action of $G$ on $\R^m$. This action allows to define algebras of vector fields which are obviously isomorphic to $so(p+1,q+1,\R)$.
As mentioned in \cite{BM}, the realization of $\g\cong
\g_{-1}\oplus\g_0\oplus\g_1$ in vector fields is the following (we denote by
$X^h$ the vector field associated to $h$) : for every $x\in\g_{-1}\cong\R^m,$
\begin{equation}\label{real}\left\{\begin{array}{ccc}
X^h_x & = & -h\quad\mbox{if}\,h\in \g_{-1}\\
X^h_x & = & -[h,x]\quad\mbox{if}\,h\in \g_{0}\\
X^h_x & = & -\frac{1}{2}[[h,x],x]\quad\mbox{if}\,h\in\g_{1}\\
\end{array}\right.\end{equation}
\subsubsection{Equivariant quantizations}
A conformally equivariant quantization (in the sense of \cite{DLO}) is a
quantization $Q :\mathcal{S}_{\delta}(\R^m)\to \mathcal{D}_{\lambda
  \mu}(\R^m)$ (where $\delta= \mu - \lambda$ is the shift value) such that, for every $X^h\in \g$, one has
\[\mathcal{L}_{X^h}\circ Q=Q\circ L_{X^h}.\]
The existence and uniqueness of such quantizations were proved in
\cite{LO,DLO} provided the shift value does not belong to a set of \emph{Critical values}.

From now on to the end of this section, we will present the tools that were
used by Duval, Lecomte and Ovsienko, and generalized in \cite{BM} in order to obtain this result.
\subsubsection{The affine quantization map}
There exists a well-known bijection from symbols to differential operators
over $\R^m$ : the so-called \emph{Standard ordering} $Q_{\mathit{Aff}}$. 
If a symbol $T\in\mathcal{S}^k_{\delta}(\R^m)$ writes
\[T(x,\xi)=\sum_{\vert\alpha\vert=k}C_{\alpha}(x)\xi^{\alpha},\]
where $\alpha$ is a multiindex, then one has
\[Q_{\mathit{Aff}}(T)=\sum_{\vert\alpha\vert=k}C_{\alpha}(x)(\frac{\partial}
{\partial x})^{\alpha}.\]
It is easily seen that this map exchanges the actions of the affine
algebra on the space of symbols
and of differential operators. This is why we call it the affine quantization 
map.

Formula \ref{real} allows to express this quantization map in a
coordinate free manner :
\begin{prop}
If $h_1,\cdots, h_k\in\R^m\cong \g_{-1}$, $A\in\Delta^{\delta}(\R^m)$, $s\in
C^{\infty}(\R^m)$, and 
\[T(x)=s(x)\,A\otimes h_1\vee\cdots\vee h_k,\]
 one has
\[Q_{\mathit{Aff}}(T) = (-1)^k s\circ A\circ L_{X^{h_1}}\circ\cdots\circ
L_{X^{h_k}},\]
where $A$ is understood as a linear map from $\Delta^{\lambda}(\R^m)$to $\Delta^{\mu}(\R^m)$.
\end{prop}
\subsubsection{The map $\gamma$}
Using the affine quantization map, one can endow the space of symbols with a
structure of representation of $\mathit{Vect}(\R^m)$ (or $\g$), 
isomorphic to
$\mathcal{D}_{\lambda \mu}$. Explicitly, we set
\[\mathcal{L}_XT:=Q_{\mathit{Aff}}^{-1}\circ\mathcal{L}_X\circ 
Q_{\mathit{Aff}}(T),\]
for every $T\in\mathcal{S}_{\delta}(\R^m)$ and $X\in\mathit{Vect}(\R^m)$.

An equivariant quantization is then a $so(p+1,q+1)$-isomorphism from
the representation $(\S, L)$ to $(\S,\L)$.

In order to measure the difference between these representations the map 
\[\gamma : \g\to gl(\S,\S) : h\mapsto \gamma(h)=\L_{X^h}-L_{X^h}\]
was introduced in \cite{BM}, where its most important properties were listed. 

In \cite{MR}, we obtained a coordinate free
expression of $\gamma$. 
\begin{prop}\label{gam}
For every $h_1,\cdots, h_k\in\R^m\cong \g_{-1}$, $A\in\Delta^{\delta}(\R^m)$
and $h\in \g_1\cong \R^{m*}$ we have
\[\begin{array}{lll}\gamma(h)(h_1\vee\cdots\vee h_k\otimes A)&=&-\lambda \sum_{i=1}^k tr([h_i,h])\,
h_1\vee\stackrel{(i)}{\cdots}\vee h_k\otimes A\\
&+&\sum_{i=1}^k\sum_{j<i}h_1\vee\stackrel{(i,j)}{\cdots}\vee
h_k\vee[h_i,[h_j,h]]\otimes A
\end{array}\]
\end{prop}
\subsubsection{Casimir operators}
In \cite{DLO,BM}, the construction of the quantization is based on the
comparison of the spectra and of the eigenvectors of some (second order)
Casimir operators. These operators are on the one hand the Casimir operator
$C$ associated to the representation $(\S, L)$ and on the other hand 
the Casimir operator $\cc$ associated to the representation
$(\S,\L)$. 

From \cite{BM}, we know that the difference between $C$ and $\cc$ can be
expressed in terms of $\gamma$ : fixing a basis $(e_i)$ in $\g_{-1}$ and
denoting by $(\epsilon^i)$ the Killing-dual basis in $\g_1$, one has
\begin{equation}\label{casinil}\cc=C+N,\end{equation}
where 
\[N=2\sum_i\gamma(\epsilon^i)L_{X^{e_i}}.\]
Moreover, it was also shown in \cite{DLO} that $\S$ is the direct sum of
eigenspaces of $C$. Indeed, if 
\[S^k(\R^m)=\oplus_{s\leq \frac{k}{2}}S_{(k,s)}\]
is the decomposition of $S^k(\R^m)$ into irreducible representations of
$\h_0\cong so(p,q)\subset \g_0$ and if we set
$\mathcal{S}_{(k,s)}=C^{\infty}(\R^m,S_{(k,s)}\otimes\Delta^{\delta}(\R^m))$,
then we have
\[C\vert_{\mathcal{S}_{(k,s)}}=\alpha_{k,s}\mbox{Id}.\]
\subsubsection{Construction of the quantization}\label{flatconst}
It turns out that the problem of existence of an $so(p+1,q+1)$-equivariant
quantization can be reduced to the following question :\\
Can we associate to every $T\in\mathcal{S}_{(k,s)}$ a \emph{unique} symbol
$\hat{T}=T_k+\cdots+T_0$ ($T_l\in \S^l(\R^m), \forall l\in\{0,\cdots,k\}$) such that
\begin{equation}\label{flatP}\left\{\begin{array}{lll}
T_k&=&T\\
\mathcal{C}(\hat{T})&=&\alpha_{k,s}\hat{T}
\end{array}\right.?
\end{equation} 
In view of (\ref{casinil}), the last equation in (\ref{flatP}) 
can be rewritten as 
\begin{equation}\label{cfin}
(C-\alpha_{k,s}\mbox{Id})T_l=N(T_{l+1})\quad\forall l\in\{0,\cdots,k-1\}
\end{equation}
In order to analyse this latter equation, the authors of \cite{DLO} introduced
the tree-like subspace associated to $\mathcal{S}_{(k,s)}$, namely
\[\widetilde{\mathcal{S}_{(k,s)}}=\bigoplus_{0\leq s-t\leq k-l}\mathcal{S}
_{(l,t)}.\]
They indeed showed that $N$ maps $\mathcal{S}_{(l,t)}$ into
$\mathcal{S}_{(l-1,t)}\oplus\mathcal{S}_{(l-1,t-1)}$. Hence, equation
(\ref{cfin}) admits a unique solution inside $\widetilde{\mathcal{S}_{(k,s)}}$ if $\alpha_{k,s}$ does not belong to the
spectrum of the restriction of $C$ to $\widetilde{\mathcal{S}_{(k,s)}}$. This
leads to the definition of critical values (see \cite{DLO}, formulas 3.7, 3.8 and 3.10 for an explicit description of the set $\Sigma_0$ of critical values).
Now the association $Q:T\mapsto \hat T$ defines an equivariant quantization because it is a bijection and fulfills
\[Q\circ L_{X^h} = \L_{X^h}\circ Q\quad\forall h\in\,\g.\]
Indeed, for all $T\in \mathcal{S}_{(k,s)}$, $Q(L_{X^h}T)$ and $ \L_{X^h}(Q(T))$ share the
following properties
\begin{itemize}
\item They are eigenvectors of $\cc$ of eigenvalue $\alpha_{k,s}$ because on
  the one hand $\cc$ commutes with $\L_{X^h}$ for all $h$ and on the other
  hand $L_{X^h}T$ belongs to $\mathcal{S}_{(k,s)}$
\item their term of degree $k$ is exactly $L_{X^h}T$,
\item they both belong to $\widetilde{\mathcal{S}_{(k,s)}}.$ 
\end{itemize}
\section{Natural and conformally equivariant quantizations up to degree 3}
We will show in this section that construction of the projectively equivariant and natural quantization of \cite{MR1} adapts in the conformal case. We will recall the sequence of propositions that leads to the existence result in \cite{MR1} and omit the proofs when the modifications are obvious. We will focus our attention to the key result that does not work for symbols of degree higher than 3.

We begin with the definition of the curved affine quantization map $Q_{\omega}$. It is based on an iterated an symmetrized version of the invariant 
differentiation 
\begin{defi}
If $f\in C^{\infty}(P,V)$ then $(\nabla^{\omega})^k f
\in C^{\infty}(P,S^k\R^{m*}\otimes V)$ is defined by
\[(\nabla^{\omega})^k f(u)(X_1,\ldots,X_k) = \frac{1}{k!}\sum_{\nu}
L_{\omega^{-1}(X_{\nu_1})}\circ\ldots\circ L_{\omega^{-1}(X_{\nu_k})}f(u)\]
for $X_1,\ldots,X_k\in\R^m$.
\end{defi}
\begin{defi}
For every symbol $T =t \otimes h_1\vee\cdots\vee h_k$, 
($t\in C^{\infty}(P,\Delta^{\delta}(\R^m))$ and $h_1,\cdots,
h_k\in\R^m\cong\g_{-1}$) we set
\begin{equation}\label{Qom}Q_{\omega}(T)f=\langle T,(\nabla^{\omega})^k f\rangle=\frac{1}{k!}\sum_{\nu}t \circ
L_{\omega^{-1}(h_{\nu_1})}\circ\cdots\circ L_{\omega^{-1}(h_{\nu_k})}f,\end{equation}
where $\nu$ runs over all permutations of the indices $\{1,\cdots,k\}$,$t$
is considered as a multiplication operator and $f$ lies in $C^{\infty}(P,\Delta^{\lambda}(\R^m))$.
\end{defi}
\begin{rem}\label{remf} If $T\in C^{\infty} 
(P,S^k_{\delta})$ is
$H-$equivariant, the differential operator $Q_{\omega}(T)$ does not transform
$H-$equivariant functions into $H-$equivariant functions; indeed, the function $Q_{\omega}(T)f$ in (\ref{Qom}) is only $G_0$-equivariant. The idea is then to  modify the
symbol $T$ by lower degree correcting terms in order to solve this problem.
\end{rem}
The following proposition shows how to measure the default of equivariance of a function defined on $P$.
\begin{prop}If $(V,\rho)$ is a representation of $G_0$ and becomes a
representation of $H$ as stated in section \ref{Lift}, then a function $v\in
 C^{\infty}(P,V)$ is $H-$equivariant iff 
\[\left\{\begin{array}{l}
\mbox{$v$ is $G_0-$equivariant}\\
\mbox{One has $L_{h^*}v=0$ for every $h$ in $\g_1$}
\end{array}\right.\]
\end{prop}
Since basically, our tools preserve the $G_0$-equivariance, we are mostly
interested in the $\g_1$-equivariance.
The following result is the keystone of our method. It works in the projective case for symbols of arbitrary degree but does not hold in the conformal case for symbols of degree higher than 3.
\begin{prop}\label{gammao3}
The relation
\begin{equation}\label{magic}L_{h^{*}}Q_{\omega}(T)(f)-Q_{\omega}(T)(L_{h^{*}}f)=Q_{\omega}
((L_{h^*}+\gamma(h))T)(f)\end{equation}
holds for all $f\in C^{\infty}(P,\Delta^{\lambda}(\R^m))_{H}$,
$h\in\mathfrak{g}_{1}$, and $T\in C^{\infty}(P,S^{k}_{\delta})$, for $k\leq 3$.
\end{prop}
\begin{proof}
It is equivalent to prove that
\[\langle T, L_{h^*}\nabla^{\omega k}f-\nabla^{\omega
  k}L_{h^*}f\rangle=\langle\gamma(h)T,\nabla^{\omega k-1}f\rangle.\]
We may check this relation on symbols of the form $T=t X^k$ where $t\in C^{\infty}(P,\Delta^\lambda(\R^m))$. Moreover, since both sides are $C^{\infty}(P)$-linear in $T$, it is sufficient to check
  this relation for a constant symbol $T$ that has the form
  $X^k$, where $X\in\mathfrak{g}_{-1}$. Then the left-hand side writes
$$L_{h*}L_{\omega^{-1}(X)}\ldots L_{\omega^{-1}(X)}f-L_{\omega^{-1}(X)}\ldots
L_{\omega^{-1}(X)}L_{h^*}f$$
and is equal to
\[\begin{array}{l}\sum_{j=1}^{k} L_{\omega^{-1}(X)}\ldots\overset{(j)}{L_{[h,X]^{*}}}\ldots L_{\omega^{-1}(X)}f\\
=\sum_{j=1}^{k-1}\sum_{i>j}L_{\omega^{-1}(X)}\overset{(j)}{\ldots}
\overset{(i)}{L_{\omega^{-1}([[h,X],X])}}\ldots L_{\omega^{-1}(X)}f\\
+\lambda \sum_{i=1}^{k}tr([h,X])L_{\omega^{-1}(X)}\overset{(i)}{\ldots}
L_{\omega^{-1}(X)}f.\end{array}\]
The first part in the latter expression is equal to
$$\frac
{k}{2}\sum_{j=1}^{k-1}L_{\omega^{-1}(X)}\ldots\stackrel{(j)}{L_{\omega^{-1}([[h,X],X])}}\ldots
L_{\omega^{-1}(X)}f$$
\begin{equation}\label{pelant}+\sum_{j=1}^{k-1}\frac {2j-k}{2} L_{\omega^{-1}(X)}\ldots
\stackrel{(j)}{L_{\omega^{-1}([[h,X],X])}}\ldots
L_{\omega^{-1}(X)}f .\end{equation}

The first sum in this expression is equal to
$$\frac {k(k-1)}{2}Q_{\omega}([[h,X],X]\vee X^{k-2}).$$
Hence, in order to obtain the desired result, we just have to show that the expression (\ref{pelant}) vanishes.

For $k$ less or equal to 2, the result is then obvious.

For $k=3$, the term (\ref{pelant}) is equal to
$$\frac {|X|^2}{2}
(L_{\omega^{-1}(X)}L_{\omega^{-1}(h^{\flat})}f-L_{\omega^{-1}(h^{\flat})}L_{\omega^{-1}(X)}f),$$
i.e. to
$$\frac {|X|^2}{2}L_{\omega^{-1}(\kappa(X,h^{\flat}))}f.$$
Using the fact that $\kappa$ has its values in
$\mathfrak{g}_{0}\oplus\mathfrak{g}_{1}$ and the $G_{1}$-equivariance
of $f$, this term is equal to
$$\frac {|X|^2}{2}L_{\om(\kappa_{0}(X,h^{\flat}))}f,$$
i.e., using the $G_{0}$-equivariance of $f$, to
$$-\frac {|X|^2}{2}\rho_{*}(\kappa_{0}(X,h^{\flat}))f.$$
This term vanishes because of the normality of $\omega$ (the trace
of $\kappa_{0}$ vanishes).

\end{proof}
\subsection{Curved Casimir operators}
We first define the analog of $N$ by setting
\[N^{\omega} : C^{\infty}(P,S^k_\delta)\to C^{\infty}(P,S^{k-1}_\delta) : T\mapsto -2\sum_{i}\gamma(\epsilon^{i})L_{\omega^{-1}(e_i)}T.\]
Then we can define the operators $C^{\omega}$ and $\mathcal{C}^{\omega}$ by
their restrictions to the spaces $C^{\infty}(P,S_{(k,s)})$ :
for all $T\in C^{\infty}(P,S_{(k,s)})$, we set
\[\left\{\begin{array}{lll}C^{\omega}(T)&=&\alpha_{k,s}T\\

\mathcal{C}^{\omega}(T)&=&C^{\omega}(T)+N^{\omega}(T),
\end{array}\right.\]
where $\alpha_{k,s}$ is the eigenvalue of $C$ on ${\mathcal S}_{(k,s)}$.

The operator $\cc^{\omega}$ has the following property
\begin{prop}\label{commute}
For every $h\in\mathfrak{g}_{1}$, one has 
$$[L_{h^{*}}+\gamma(h), \mathcal{C}^{\omega}]=0$$
on $C^{\infty}(P,S^k_\delta)_{G_0}$.
\end{prop}
For the operator $N^{\omega}$, we have the following result :
\begin{prop}\label{Goinv}
The operator $N^{\omega}$ preserves the $G_0$-equivariance of functions.
\end{prop}
\section{Construction of the quantization}\label{curvconst}
First remark that the construction of section \ref{flatconst} is still valid
in the curved case.
\begin{thm}\label{hat}
If $\delta$ does not belong to the set $\Sigma_0$ of critical values , for every $T$
in $C^{\infty}(P,S_{(k,s)})$, there exists a unique function $\hat{T}$ in $ C^{\infty}(P,\widetilde{S_{(k,s)}})$
 such that
\begin{equation}\label{curvP}\left\{\begin{array}{lll}
\hat{T}&=&T_k+\cdots+T_0,\quad T_k=T\\
\mathcal{C}^{\omega}(\hat{T})&=&\alpha_{k,s}\hat{T}.
\end{array}\right.
\end{equation} 
Moreover, if $T$ is $G_0$-invariant, then $\hat{T}$ is $G_0$-invariant.
\end{thm}
This result allows to define the main ingredient in order to define the
quantization.
\begin{defi}\label{Qflat}
Suppose that $\delta$ is not critical. Then
the map 
\[Q : C^{\infty}(P,S_{\delta})\to  C^{\infty}(P,S_{\delta})\]
is the linear extension of the association $T\mapsto \hat{T}$.
\end{defi}
The map $Q$ has the following  property :
\begin{prop}\label{propQ}
There holds
\begin{equation}
\label{Q}(L_{h^*}+\gamma(h))Q(T)=Q(L_{h^*}T),
\end{equation}
for every $h\in\g_1$ and every $T\in C^{\infty}(P,S_{\delta})_{G_0}$.
\end{prop}
Finally, we obtain
\begin{thm}
If $\delta$ does not belong to the set $\Sigma_0$ of critical values (\cite[Formula 3.10]{DLO}), then the formula
\begin{equation}\label{quant}Q_M: (g,T)\mapsto Q_M(g,T)(f)=(p^*)^{-1}[Q_{\omega}(Q(p^*T))(p^*f)],\end{equation}
(where $Q_{\omega}$ is given by (\ref{Qom})) defines a natural and conformally equivariant quantization on symbols of degree at most three.
\end{thm}
\section{Quantization of symbols of degree four}
In this section, we will first show how proposition \ref{gammao3} has to be adapted for symbols of degree four in the presence of a nontrivial curvature. We will then introduce the corresponding correcting terms to the map $Q$ of definition \ref{Qflat} in order to prove the existence of the natural and conformally equivariant quantization for symbols of degree four.
\subsection{Modification of the operator $\gamma$}
The key result in the construction of the quantization is proposition \ref{gammao3}. For symbols of degree four we introduce the following operators.
\begin{defi} We define the operators
\[\gamma_3 :\g_1\otimes C^{\infty}(P,S^4_\delta)\to C^{\infty}(P,S^1_\delta)\]
and 
\[\gamma_4 :\g_1\otimes C^{\infty}(P,S^4_\delta)\to C^{\infty}(P,S^0_\delta)\]
by letting $\gamma_3(h)$ and $\gamma_4(h)$ be $C^{\infty}(P)$-linear for every $h$ in $\g_1$ and by defining the values of these operators on a symbol of the form $X^4$
\[\left\{\begin{array}{lll}
\gamma_3(h)(X^4)&=&|X|^2\kappa_0(h^\flat,X)X\\
\gamma_4(h)(X^4)&=&-\lambda m \,|X|^2\langle\kappa_1(h^\flat,X),X\rangle
\end{array}\right.\]
On the spaces of symbols of degree less or equal to 3, we also set $\gamma_3(h)=\gamma_4(h)=0$.
\end{defi}
Finally
we define the modified operator $\gamma'$ on symbols of degree less or equal to four :
\begin{defi}
The map $\gamma'$ is defined on $\oplus_{k=0}^4\g_1\otimes C^{\infty}(P,S^k_\delta)$ by
\[\gamma'(h)=\gamma(h)+\gamma_3(h)+\gamma_4(h),\]
for every $h$ in $\g_1$.
\end{defi}
We then have the following result, which is the adaptation of Proposition \ref{gammao3}.
\begin{prop}\label{gammao4}
There holds 
\begin{equation}\label{magic4}L_{h^{*}}Q_{\omega}(T)(f)-Q_{\omega}(T)(L_{h^{*}}f)=Q_{\omega}
((L_{h^*}+\gamma'(h))T)(f)\end{equation}
 for all $f\in C^{\infty}(P,\Delta^{\lambda}(\R^m))_{H}$,
$h\in\mathfrak{g}_{1}$, and $T\in C^{\infty}(P,S^{k}_{\delta})$, for $k\leq 4$.
\end{prop}
\begin{proof}
The first part of the proof of proposition \ref{gammao3} is still valid for symbols of degree 4. It now remains to compute the expression (\ref{pelant}). It is equal to
\[|X|^2\left[L_{\om(X)}L_{\om(X)}L_{\om(h^\flat)}f-L_{\om(h^\flat)}L_{\om(X)}L_{\om(X)}f\right]. \]
The next step is to use Lemma \ref{omega} to compute the second term. We get that expression (\ref{pelant}) is equal to
\[|X|^2\left[L_{\om(X)}L_{\om(\kappa(h^\flat,X))}f+L_{\om(\kappa(h^\flat,X))}L_{\om(X)}f\right].\]
The first term of this expression vanishes due to the $H$-invariance of $f$ (recall that $\kappa$ has values in $\h$). The second term is easy to compute using Lemma \ref{omega} and the invariance of $f$.
\end{proof}
\subsection{Modification of the application $Q$}
In the construction of the quantization up to degree 3, the main interest of the application $Q$ (see Definition \ref{Qflat}) is that the operator $Q_{\omega}(Q(T))$, for all $H$-equivariant functions $T$, transforms $H$-equivariant functions into $H$-equivariant functions. This ensures that formula (\ref{quant}) makes sense. The following proposition shows how this property fails for symbols of degree 4.
\begin{prop}\label{derQ}
There holds 
\[L_{h^*}[(Q_{\omega}(Q(T)))(f)]= Q_{\omega}(\gamma_3(h)T+\gamma_4(h)T)(f)\]
for all $h\in\g_1$, $T\in C^{\infty}(P,S^k_{\delta})_H$ ($k\leq 4$) and $f\in C^{\infty}(P,\Delta^\lambda(\R^m))_H$.
\end{prop}
\begin{proof}
By Proposition \ref{gammao4}, we have 
\begin{equation}\label{failQ}L_{h^*}[(Q_{\omega}(Q(T)))(f)]=Q_{\omega}(Q(T))(L_{h^*}f)+Q_{\omega}((L_{h^*}+\gamma'(h))Q(T))(f).\end{equation}
The first term is vanishing. In view of Proposition \ref{propQ}, the second one is equal to
\[Q_{\omega}(Q(L_{h^*}T))+Q_{\omega}((\gamma_3(h)+\gamma_4(h))Q(T)).\]
The result follows from the $H$-invariance of $T$ and from the vanishing of $\gamma_3(h)$ and $\gamma_4(h)$ on symbols of degree less or equal to 3.
\end{proof}
We will now introduce some modifications to the map $Q$ in order to annihilate the right hand side of equation (\ref{failQ}). In the next definition, we still denote by $(e_j)$ a basis of $\g_{-1}\cong \R^m$ and we denote by $(\eta^j)$ the dual basis in the usual sense.

Recall that the divergence associated to a Cartan connection is then defined by
\[div^\omega : C^{\infty}(P,S^k_\delta)\to C^{\infty}(P,S^{k-1}_\delta) : T\mapsto \sum_j L_{\om(e_j)}i(\eta^j)T.\]
\begin{defi}When $\delta\not\in\{\frac{m+1}{m},\frac{m+2}{m}\}$, we define the maps $Q_3$ and $Q_4$ explicitly on symbols of the form $tX^4$ where $t\in C^{\infty}(P,\Delta^{\delta}(\R^m))$ and $X\in\g_{-1}$ by
\[Q_3(tX^4)=-|X|^2\left[t\langle\kappa_1(\eta^{j^\flat},X),X\rangle e_j+\frac{2}{m+2-m\delta}
div^{\omega}(t\kappa_0(\eta^{j^\flat},X)X\vee e_j)\right]\]
\[Q_4(tX^4)=\frac{-m\lambda}{(m+1-m\delta)(m+2-m\delta)}|X|^2div^{\omega^2}(t\kappa_0(\eta^{j^\flat},X)X\vee e_j).\vspace{0.2cm}\]
We also set $Q'=Q+Q_3+Q_4$.
\end{defi}
\noindent {\bf Remark :} It is easy to check that the values $\frac{m+1}{m}$ and $\frac{m+2}{m}$ belong to the set $\Sigma_0$ of critical values, so that the formulas for $Q_3$ and $Q_4$ make sense if $\delta$ is not critical.

In order to obtain the analog of Proposition \ref{propQ}, we need the following technical lemma.
\begin{lem}\label{div}
For every $T\in C^{\infty}(P,S^k_\delta)_{G_0}$ and $h\in\g_1$, one has
\[L_{h^*}div^{\omega}T-div^\omega L_{h^*}T=(m+2k-2-m\delta)i(h)T-(k-1)i(\eta^j)i(e_j^\sharp)T\vee h^\flat.\]
For $T\in C^{\infty}(P,S^2_\delta)_{G_0}$ and $h\in\g_1$, one has
\[L_{h^*}div^{\omega^2}T-div^{\omega^2} L_{h^*}T=2(m+1-m\delta)i(h)div^{\omega}T-L_{\om(h^\flat)}i(\eta^j)i(e_j^\sharp)T.\]
\end{lem}
\begin{proof}
We proceed as in \cite[Lemma 8]{MR} : we have for $\xi^1,\cdots,\xi^{k-1}\in\R^{m^*}$
\[\begin{array}{l}L_{h^*}div^{\omega}T(\xi^1,\cdots,\xi^{k-1})-div^\omega L_{h^*}T(\xi^1,\cdots,\xi^{k-1})=\\(\rho(h\otimes e_j+h_j Id-e_j^\sharp\otimes h^\flat)T)(\eta^j,\xi^1,\cdots,\xi^{k-1})=\\
\left[(m+2k-1-(m+1)\delta)+(\delta-1)\right]i(h)T(\xi^1,\cdots,\xi^{k-1})\\ -\sum_l\langle \xi^l,h^\flat\rangle T(\eta^j,e_j^\sharp,\xi^1\stackrel{(l)}{\cdots},\eta^{k-1}).\end{array}\]
This yields the first part of the result. The second part follows by induction.
\end{proof}
\begin{prop}\label{prop20}
There holds
\[(L_{h^*}+\gamma(h))(Q_3+Q_4)(T)=-(\gamma_3(h)+\gamma_4(h))T\]
for every $h\in \g_1$ and $T\in C^{\infty}(P,S^4_\delta)_H$.
\end{prop}
\begin{proof}
The contributions of the term containing $\gamma(h)$ are easy to compute because $\gamma(h)$ vanishes on symbols of degree 0 and reduces to $-\lambda m \,i(h)$ on symbols of degree 1.

We successively use Lemmas \ref{div}, \ref{derkappa} and the normality of the Cartan connection to compute the terms in $L_{h^*}$ and get the desired result. 
\end{proof}
We then have an important corollary.
\begin{cor}\label{cor}
For every $f\in C^{\infty}(P,\Delta^{\lambda}(\R^m))_H$ and $T\in C^{\infty}(P,S^4_\delta)_H$, the function
$Q_{\omega}(Q'(T))(f)$ is $H$-invariant.
\end{cor}
\begin{proof}
First notice that the map $Q'$ transforms $G_0$-invariant functions into $G_0$-invariant functions. Indeed, this is certainly true for $Q$ as quoted in \cite{MR1}, and holds true for $Q_3$ and $Q_4$, because of the properties of invariance of $\kappa_0$ and $\kappa_1$ (see \cite{Capinv}). Therefore, the operator $Q_{\omega}(Q'(T))$ also transforms $G_0$-equivariant functions into $G_0$ equivariant functions. We just need to show that $Q_{\omega}(Q'(T))(f)$ is $\g_1$-invariant. This follows immediately from Propositions \ref{derQ} and \ref{prop20}.
\end{proof}
Finally, we have the main result :
\begin{thm}
If $\delta$ is not critical, then the formula
\begin{equation}Q_M: (g,T)\mapsto Q_M(g,T)(f)=(p^*)^{-1}[Q_{\omega}(Q'(p^*T))(p^*f)],\end{equation}
defines a natural and conformally equivariant quantization on symbols of degree four\end{thm}
\begin{proof}
The proof goes as in \cite{MR1}. The main point is that the formula makes sense. But this is a consequence of corollary \ref{cor}.
\end{proof}
{\bf Acknowledgements :} F. Radoux is supported by BFR 06/077 from the Minist\`{e}re de la Culture, de l' Enseignement sup\'erieur et de la Recherche of the Grand Duchy of Luxembourg.
\bibliographystyle{plain} \bibliography{prquant}

\begin{thebibliography}{10}

\bibitem{BHMP}
F.~Boniver, S.~Hansoul, P.~Mathonet, and N.~Poncin.
\newblock Equivariant symbol calculus for differential operators acting on
  forms.
\newblock {\em Lett. Math. Phys.}, 62(3):219--232, 2002.

\bibitem{BM}
F.~Boniver and P.~Mathonet.
\newblock Ifft-equivariant quantizations.
\newblock {\em To appear in J. Geom. Phys., math.RT/0206213}, 2005.

\bibitem{Bor}
M.~Bordemann.
\newblock Sur l'existence d'une prescription d'ordre naturelle projectivement
  invariante.
\newblock {\em Submitted for publication, math.DG/0208171}.

\bibitem{Bou1}
Sofiane Bouarroudj.
\newblock Projectively equivariant quantization map.
\newblock {\em Lett. Math. Phys.}, 51(4):265--274, 2000.

\bibitem{Bou2}
Sofiane Bouarroudj.
\newblock Formula for the projectively invariant quantization on degree three.
\newblock {\em C. R. Acad. Sci. Paris S\'er. I Math.}, 333(4):343--346, 2001.

\bibitem{Capinv}
A.~{\v{C}}ap, J.~Slov{\'a}k, and V.~Sou{\v{c}}ek.
\newblock Invariant operators on manifolds with almost {H}ermitian symmetric
  structures. {I}. {I}nvariant differentiation.
\newblock {\em Acta Math. Univ. Comenian. (N.S.)}, 66(1):33--69, 1997.

\bibitem{Capcart}
A.~{\v{C}}ap, J.~Slov{\'a}k, and V.~Sou{\v{c}}ek.
\newblock Invariant operators on manifolds with almost {H}ermitian symmetric
  structures. {II}. {N}ormal {C}artan connections.
\newblock {\em Acta Math. Univ. Comenian. (N.S.)}, 66(2):203--220, 1997.

\bibitem{DLO}
C.~Duval, P.~Lecomte, and V.~Ovsienko.
\newblock Conformally equivariant quantization: existence and uniqueness.
\newblock {\em Ann. Inst. Fourier (Grenoble)}, 49(6):1999--2029, 1999.

\bibitem{DO1}
C.~Duval and V.~Ovsienko.
\newblock Conformally equivariant quantum {H}amiltonians.
\newblock {\em Selecta Math. (N.S.)}, 7(3):291--320, 2001.

\bibitem{DO}
C.~Duval and V.~Ovsienko.
\newblock Projectively equivariant quantization and symbol calculus:
  noncommutative hypergeometric functions.
\newblock {\em Lett. Math. Phys.}, 57(1):61--67, 2001.

\bibitem{sarah}
Sarah Hansoul.
\newblock Existence of natural and projectively equivariant quantizations.
\newblock {\em To appear, math.DG/0601518.}

\bibitem{Hansoul}
Sarah Hansoul.
\newblock Projectively equivariant quantization for differential operators
  acting on forms.
\newblock {\em Lett. Math. Phys.}, 70(2):141--153, 2004.

\bibitem{Kobabook}
Shoshichi Kobayashi.
\newblock {\em Transformation groups in differential geometry}.
\newblock Springer-Verlag, New York, 1972.
\newblock Ergebnisse der Mathematik und ihrer Grenzgebiete, Band 70.

\bibitem{LO}
P.~B.~A. Lecomte and V.~Yu. Ovsienko.
\newblock Projectively equivariant symbol calculus.
\newblock {\em Lett. Math. Phys.}, 49(3):173--196, 1999.

\bibitem{Leconj}
Pierre B.~A. Lecomte.
\newblock Towards projectively equivariant quantization.
\newblock {\em Progr. Theoret. Phys. Suppl.}, (144):125--132, 2001.
\newblock Noncommutative geometry and string theory (Yokohama, 2001).

\bibitem{Loubon}
S.~E. Loubon~Djounga.
\newblock Conformally invariant quantization at order three.
\newblock {\em Lett. Math. Phys.}, 64(3):203--212, 2003.

\bibitem{MR}
Pierre Mathonet and Fabian Radoux.
\newblock Natural and projectively equivariant quantiations by means of cartan
  connections.
\newblock {\em Lett. Math. Phys.}, 72(3):183--196, 2005.

\bibitem{MR1}
Pierre Mathonet and Fabian Radoux.
\newblock Cartan connections and natural and projectively equivariant
  quantizations.
\newblock {\em J. London Math. Soc.}, To appear, 2007.

\bibitem{Woo}
N.~M.~J. Woodhouse.
\newblock {\em Geometric quantization}.
\newblock Oxford Mathematical Monographs. The Clarendon Press Oxford University
  Press, New York, second edition, 1992.
\newblock Oxford Science Publications.

\end{thebibliography}
\end{document}